\documentclass{article}\begin{document}
\centerline{\bf On Generalized Integrals and Ramanujan-Jacobi Special Functions}
\vskip .4in   

\centerline{N.D. Bagis}
\centerline{Department of Informatics, Aristotele University}
\centerline{Thessaloniki, Greece}
\centerline{nikosbagis@hotmail.gr}
\vskip .2in

\[
\]

\textbf{Keywords}: Integrals; Elliptic Functions; Ramanujan; Special Functions; Continued Fractions; Generalization; Evaluations
\[
\]
\centerline{\bf Abstract}

\begin{quote}
In this article we consider new generalized functions for evaluating integrals and roots of functions. The construction of these generalized functions is based on Rogers-Ramanujan continued fraction, the Ramanujan-Dedekind eta, the elliptic singular modulus and other similar functions. We also provide modular equations of these new generalized functions and remark some interesting properties.   

\end{quote}

\section{Introduction}

Let 
\begin{equation}
\eta(\tau)=e^{\pi i\tau/12} \prod^{\infty}_{n=1}(1-e^{2\pi in\tau})
\end{equation}
denotes the Dedekind eta function which is defined in the upper half complex plane. It not defined for real $\tau$.\\
Let for $|q|<1$ the Ramanujan eta function be 
\begin{equation}
f(-q)=\prod^{\infty}_{n=1}(1-q^n).
\end{equation}
The following evaluation holds (see [12]):
\begin{equation}
f(-q)=2^{1/3}\pi^{-1/2}q^{-1/24}k^{1/12}k^{{*}{1/3}}K(k)^{1/2},
\end{equation}
where $k=k_r$ is the elliptic singular modulus, $k^{*}=\sqrt{1-k^2}$ and $K(x)$ is the complete elliptic integral of the first kind.\\
The Rogers-Ramanujan continued fraction is (see [8],[9],[14]):
\begin{equation}
R(q):=\frac{q^{1/5}}{1+}\frac{q^1}{1+}\frac{q^2}{1+}\frac{q^3}{1+}...  ,
\end{equation}
which have first derivative (see [5]):
\begin{equation}
R'(q)=5^{-1}q^{-5/6}f(-q)^4R(q)\sqrt[6]{R(q)^{-5}-11-R(q)^5}
\end{equation}
and also we can write 
\begin{equation}
\frac{dR(q)}{dk}=5^{-1}\cdot2^{1/3}(kk^{*})^{-2/3}R(q)\sqrt[6]{R(q)^{-5}-11-R(q)^5}.
\end{equation}
Ramanujan have proven that
\begin{equation}
R(q)^{-5}-11-R(q)^5=\frac{f(-q)^6}{qf(-q^5)^6}.
\end{equation} 
Also holds the following interesting identity
\begin{equation}
\frac{dk}{dq}=\frac{-2k(k^{*})^2K(k)^2}{q\pi^2},
\end{equation}   
which is a result of Ramanujan (see [7] and [8]) for the first derivative of $k=k_r$ with respect to $q=e^{-\pi\sqrt{r}}$, $r>0$. 

\section{Propositions}

\textbf{Definition 1.}\\
For any smooth function $G$, we define $m_G(A)$ to be such that
\begin{equation}
A=\pi\int^{+\infty}_{\sqrt{m_G(A)}}\eta\left(it/2\right)^4G\left(R\left(e^{-\pi t}\right)\right)dt.
\end{equation}
\\
\textbf{Theorem 1.}\\
If 
\begin{equation}
5\int^{y}_{0}\frac{G(t)}{t\sqrt[6]{t^{-5}-11-t^5}}dt=A,
\end{equation}
then
\begin{equation}
y=R\left(e^{-\pi\sqrt{m_G(A)}}\right).
\end{equation}
\\
\textbf{Proof.}\\
From [5] it is known that if $a,b\in (0,1)$ then
\begin{equation}
\int^{a}_{b}f(-q)^4q^{-5/6}G(R(q))dq=5\int^{R(a)}_{R(b)}\frac{G(x)}{x\sqrt[6]{x^{-5}-11-x^5}}dx ,
\end{equation}
which is equivalent to write
$$
-\pi \int^{a_1}_{b_1}f\left(-e^{-t\pi}\right)^4e^{-t\pi/6}G\left(R\left(
e^{-\pi t}\right)\right)dt=5\int^{R\left(e^{-\pi a_1}\right)}_{R\left(e^{-\pi b_1}\right)}\frac{G(x)}{x\sqrt[6]{x^{-5}-11-x^5}}dx .
$$
Setting $b_1=+\infty$, $a_1=\sqrt{m_G(A)}$ and then using Definition 1 and the Dedekind eta expansion (1) we get the result.\\

Also differentiating (8) one can get
\begin{equation}
\frac{dm_G^{(-1)}(r)}{dr}=-\frac{\pi\eta\left(i\sqrt{r}/2\right)^4G(R(q))}{2\sqrt{r}}\textrm{, }q=e^{-\pi\sqrt{r}}\textrm{, }r>0
\end{equation}
and if 
\begin{equation}
\phi(r):=-\frac{2\sqrt{r}}{\pi \eta\left(\frac{i\sqrt{r}}{2}\right)^4G\left(R\left(q\right)\right)},
\end{equation}
then\\
\\
\textbf{Proposition 1.}\\
If $A>0$, then
\begin{equation}
\frac{d}{dA}m_G(A)=\phi\left(m_G(A)\right).
\end{equation}
\\

Derivating (10) with respect to $q_{m_G}=e^{-\pi\sqrt{m_G}}$ we get
$$
\frac{5G(R(q_{m_G}))}{R(q_{m_G})\sqrt[6]{R(q_{m_G})^{-5}-11-R(q_{m_G})^5}}\frac{dR(q_{m_G})}{dq_{m_G}}=\frac{dA}{dk_{m_G}}\frac{1}{\frac{dq_{m_G}}{dk_{m_G}}},
$$
or the equivalent, using (8),(5) and (3):\\
\\
\textbf{Theorem 2.}
\begin{equation}
\frac{dA}{dk_{m_G}}=\sqrt[3]{2}\frac{G\left(y(A)\right)}{\left(k_{m_G}k^{*}_{m_G}\right)^{2/3}}.
\end{equation}
\\

Also with integration of (16) we have
\begin{equation}
3\sqrt[3]{2k_{m_G}}\cdot{}_2F_{1}\left[\frac{1}{6},\frac{1}{3};\frac{7}{6};k_{m_G}^2\right]=\int^{A}_{0}\frac{1}{G\left( R\left(e^{-\pi\sqrt{m_G(t)}}\right)\right)}dt.
\end{equation}
Hence\\
\\
\textbf{Definition 2.}\\
We define the function $h$ as
\begin{equation}
x=\int^{h(x)}_{0}\frac{dt}{G\left(R\left(e^{-\pi\sqrt{m_G(t)}}\right)\right)}=\int^{h(x)}_{0}\frac{dt}{G\left(y(t)\right)}
\end{equation}
\textbf{Theorem 3.}\\
Set $m_G(A)=r$, then
\begin{equation}
A=h\left(3\sqrt[3]{2k_r}\cdot{}_2F_{1}\left[\frac{1}{6},\frac{1}{3};\frac{7}{6};k_r^2\right]\right),
\end{equation}
or beter
\begin{equation}
m_G^{(-1)}(r)=h\left(3\sqrt[3]{2k_r}\cdot{}_2F_{1}\left[\frac{1}{6},\frac{1}{3};\frac{7}{6};k_r^2\right]\right).
\end{equation}
\\

From relation (10) differentiating we have
$$
5\frac{G(y(A))y'(A)}{y(A)\sqrt[6]{y(A)^{-5}-11-y(A)^5}}=1.
$$
Hence
\begin{equation}
5\int^{y(A)}_{0}\frac{dt}{t\sqrt[6]{t^{-5}-11-t^5}}=\int^{A}_{0}\frac{dt}{G(y(t))}=b_{m_G},
\end{equation}
where 
\begin{equation}
b_r=3\sqrt[3]{2k_r}\cdot{}_2F_{1}\left[\frac{1}{6},\frac{1}{3};\frac{7}{6};k_r^2\right].
\end{equation}
Also if $m(x)$ is the function defined as (see [13])
\begin{equation}
\pi\int^{+\infty}_{\sqrt{m(r)}}\eta(it/2)^4dt=r.
\end{equation}
Then we have also 
\begin{equation}
b^{(-1)}_r=m(r).
\end{equation}
Hence
\begin{equation}
y(A)=R\left(e^{-\pi\sqrt{m_G}}\right)=F_1\left(b_{m_G}\right),
\end{equation}
where we have define $F_1$ by
\begin{equation}
x=5\int^{F_1(x)}_{0}\frac{dt}{t\sqrt[6]{t^{-5}-11-t^5}}.
\end{equation}
\textbf{Theorem 4.}
\begin{equation}
5\int^A_0\frac{G(t)dt}{t\sqrt[6]{t^{-5}-11-t^5}}=h\left(5\int^{A}_{0}\frac{dt}{t\sqrt[6]{t^{-5}-11-t^5}}\right).
\end{equation}
\textbf{Proof.}\\
From relations (21) and (19) and (10) we have
\begin{equation}
h\left(5\int^{y(A)}_{0}\frac{dt}{t\sqrt[6]{t^{-5}-11-t^5}}\right)=h(b_{m_{G}})=A=5\int^{y(A)}_{0}\frac{G(t)dt}{t\sqrt[6]{t^{-5}-11-t^5}},
\end{equation}
which by inversion of $y(A)$ we get the desired result.\\

Also from relation (20) is
\begin{equation}
m^{(-1)}_G(A)=h(b_A).
\end{equation}
Hence knowing the function $h$ we know almost everything. For this, for a given function $P(x)$ we can set 
\begin{equation}
G(x)=5^{-1}xP'(x)\sqrt[6]{x^{-5}-11-x^5},
\end{equation} 
then 
\begin{equation}
m_G^{(-1)}(A)=P\left(R\left(e^{-\pi\sqrt{A}}\right)\right)=h(b_A).
\end{equation}
Hence inverting $b_A$ we get the function $h$. Note also that $P(A)=y^{(-1)}(A)$ and
\begin{equation}
m_G(A)=v_i(y(A)),
\end{equation}
where $v_i(x)$ is the inverse function of $R(e^{-\pi\sqrt{x}})$.\\
\\
\textbf{Remarks.}\\ 
1) Note by definition that when we know $G$ the $m_G^{(-1)}(A)$ is a closed form formula and $m_G(A)$ is not  (it needs inversion).\\
2) Continuing from relation (30) and $y(x)=P^{(-1)}(x)$ we have
$$
G(y(x))=5^{-1}P^{(-1)}(x)P'\left(P^{(-1)}(x)\right)\sqrt[6]{\left(P^{(-1)}(x)\right)^{-5}-11-\left(P^{(-1)}(x)\right)^5}.
$$
From equation (28) we get also\\
\\
\textbf{Corollary 1.}\\
\begin{equation}
h^{(-1)}(A)=5\int^{y(A)}_{0}\frac{dt}{t\sqrt[6]{t^{-5}-11-t^5}}.
\end{equation} 
\\
\textbf{Corollary 2.}
\begin{equation}
y'(A)=5^{-1}h'_i(A)y(A)\sqrt[6]{y(A)^{-5}-11-y(A)^5},
\end{equation}
with 
\begin{equation}
h_i(A)=h^{(-1)}(A)=b_{m_G(A)}=\left(b\circ m_G\right)(A)
\end{equation}
and also
\begin{equation}
\frac{dy(A)}{dk_{m_G}}=5^{-1}\cdot 2^{1/3}\left(k_{m_G}k^{*}_{m_G}\right)^{-2/3}y(A)\sqrt[6]{y(A)^{-5}-11-y(A)^5}.
\end{equation}
\\
Using the above one can show with differentiation the following\\
\\
\textbf{Theorem 5.}
\begin{equation}
h_i(A)=\pi\int^{\infty}_{\sqrt{m_G(A)}}\eta\left(it/2\right)^4dt=b_{m_G(A)}=5\int^{y(A)}_{0}\frac{dt}{t\sqrt[6]{t^{-5}-11-t^5}}.
\end{equation}
\\

Also another interesting theorem arises from the definition of $h(x)$. By setting
\begin{equation}
h_1(t):=\left(\frac{1}{h_i'(x)}\right)^{(-1)}(t),
\end{equation}
then\\
\\
\textbf{Theorem 6.}
\begin{equation}
5\int^{G_i(x)}_{0}\frac{dt}{t\sqrt[6]{t^{-5}-11-t^5}}=\int^{x}_{c}\frac{h'_1(t)}{t}dt.
\end{equation}
\\
\textbf{Proof.}\\
From relation (18) we have
$$
G\left(y(t)\right)=\frac{1}{h_i'(t)},
$$
hence inverting
$$
y_i\left(G_i(t)\right)=h_1(t).
$$
Taking the derivatives in both parts we get
$$
5\frac{G_i'(t)}{G_i(t)\sqrt[6]{G_i(t)^{-5}-11-G_i(t)^5}}=\frac{h_1'(t)}{t}.
$$
Lastly integrating the above relation in both parts we get the result.\\

Equation (39) can also written in the next form using (26)
\begin{equation}
G_i(x)=F_1\left(\int^{x}_{c}\frac{h_1'(t)}{t}dt\right).
\end{equation}
  
Suppose the $n$-th order modular equation of $y(A)$ is
\begin{equation}
P_n(A)=y\left(n\cdot y^{(-1)}(A)\right).
\end{equation}
Then we can find $P_n(x)$ by solving
\begin{equation}
\int^{P_n(A)}_{0}\frac{G(t)}{t\sqrt[6]{t^{-5}-11-t^5}}dt=n\int^{A}_{0}\frac{G(t)}{t\sqrt[6]{t^{-5}-11-t^5}}dt,
\end{equation}
with respect to $P_n(A)$. Setting $Q_n(A)$ to be 
\begin{equation}
Q_n(A):=m_G\left(n^2\cdot  m_G^{(-1)}\left(A\right)\right),
\end{equation}
the $n$-nth degree modular equation of $m_G(A)$ and using (15) we have 
\begin{equation}
m_G^{(-1)}(A)=\int^{A}_{0}\frac{dt}{\phi(t)}
\end{equation}
and 
$$
\int^{Q_n(A)}_{0}\frac{dt}{\phi(t)}=n\int^{A}_{0}\frac{dt}{\phi(t)} .
$$
Also
$$
Q^{(-1)}_n\left(A\right)=Q_{1/n}\left(A\right) 
$$
and
$$
\int^{Q_{nm}(A)}_{0}\frac{dt}{\phi(t)}=n\int^{Q_m(A)}_{0}\frac{dt}{\phi(t)}=nm\int^{A}_{0}\frac{dt}{\phi(t)}
$$
and
$$
\int^{Q_{n}\left(Q_{m}\left(A\right)\right)}_{0}\frac{dt}{\phi(t)}=n\int^{Q_m(A)}_{0}\frac{dt}{\phi(t)}=nm\int^{A}_{0}\frac{dt}{\phi(t)}.
$$
Hence
$$
\int^{Q_{n}\left(Q_{m}\left(A\right)\right)}_{0}\frac{dt}{\phi(t)}=\int^{Q_{nm}(A)}_{0}\frac{dt}{\phi(t)}
$$
and consequently
\begin{equation}
Q_n\left(Q_m\left(A\right)\right)=Q_m\left(Q_n\left(A\right)\right)=Q_{nm}(A),
\end{equation}
$Q_1(A)=A$
and if $n=p_1^{a_1}p_2^{a_2}\ldots p_s^{a_s}$
\begin{equation}
Q_n(x)=\left(Q_{p_1}\circ\ldots\circ Q_{p_1}\right)\circ \left(Q_{p_2}\circ \ldots \circ Q_{p_2}\right)\circ\ldots \circ\left( Q_{p_s}\circ\ldots \circ Q_{p_s}\right),
\end{equation}
respectively iterated $a_1,a_2,\ldots,a_s$ times.\\ 
Hence
$$
y\left(m^{(-1)}_G\circ m_G^{(-1)}\circ Q_{n}\circ m_G\circ m_G(A)\right)=R\left(e^{-\pi\sqrt{m^{(-1)}_G\circ Q_{n}\circ m_G\circ m_G(A)}}\right)=
$$
$$
=R\left(e^{-\pi n\sqrt{m_G(A)}}\right)=\Omega_n\left(y(A)\right).
$$
But 
\begin{equation}
m^{(-1)}_G\circ m_G^{(-1)}\circ Q_{n}\circ m_G\circ m_G(A)=m_G^{(-1)}\left(n^2\cdot m_G(A)\right).
=Q^{*}_{n}(A)
\end{equation}
\\
By this we lead to the following\\
\\
\textbf{Theorem 7.}
$$
y\left(Q^{*}_{n^2}(A)\right)=\Omega_n\left(y(A)\right),
$$
where $\Omega_n(A)$ is the $n$-th modular equation of the Rogers-Ramanujan continued fraction and $Q_{n^2}^{*}(A)$ is the $n$-th order modular equation of $m_G^{(-1)}(A)$ i.e $Q^{*}_{n^2}(A)=m_G^{(-1)}\left(n^2\cdot  m_G\left(A\right)\right)$.\\ 
\\
\textbf{Theorem 8.}\\If we know $f_0=y$ and $y_i=f^{(-1)}_0$, then 
$$
m_G(A)=b^{(-1)}\circ F_1^{(-1)}\circ f_0(A)\eqno{:(\nu1)}
$$
and for all $n$ positive real: 
$$
m_G^{(-1)}(n^2)=f^{(-1)}_0\circ\Omega_{n}\circ f_0\circ m_G^{(-1)}(1)=f^{(-1)}_0\left(R\left(e^{-\pi n}\right)\right),\eqno{:(\nu 2)}
$$
where $n$ can take and positive real values as long as
$$
\Omega_{n}(x)=v\left(n^2\cdot v^{(-1)}(x)\right)\eqno{:(\nu 3)}
$$
and $v(x)=R\left(e^{-\pi\sqrt{x}}\right)$, i.e. $v(x)$ is the Rogers-Ramanujan continued fraction. The constant function $\Omega_n(x)$ is the $n$-th degree modular equation of the $v(x)$. Hence knowing $y(x)$ and $y^{(-1)}(x)=y_i(x)$ we know respectively $m_G(x)$ and $m_G^{(-1)}(x)$.\\ 
\\
\textbf{Examples.}\\
\textbf{1)} If $f(x)=x^2+2x$ then $f^{(-1)}(x)=-1+\sqrt{1+x}$ and 
$$
m_G(x)=b^{(-1)}\circ F_1^{(-1)}\left(-1+\sqrt{1+x}\right).
$$
\\
\textbf{2)} If 
$$
G(t)=\frac{t\sqrt[6]{t^{-5}-11-t^5}}{t+1},
$$ 
then $y(x)=-1+e^{x/5}$ and $y_i(x)=5\log(x+1)$. Hence exist constant $c_1$ such that
$$
m_G^{(-1)}(n^2)=5\log\left(1+\Omega_{n}(c_1)\right)\textrm{, }\forall n>0
$$
\\
\textbf{3)} 
If 
$$
G(x)=\frac{x\sqrt[6]{x^{-5}-11-x^5}}{5\sqrt{1-k \sin(x)^2}},
$$
then
$$
y(x)=E(x,k)\textrm{ and }y_i(x)=\textrm{am}(x,k),
$$
where $\textrm{am}$ is the Jacobi amplitude i.e the inverse of the incomplete elliptic integral of the first kind $E[x,k]$, with $k$ a parameter. Hence
$$
m_G^{(-1)}(4x)=Q^{*}_4\left(m_G^{(-1)}(x)\right)\textrm{, }Q^{*}_4(x)=E\left[\Omega_2\left(\textrm{am}(x,k)\right),k\right].
$$
If $k=1$
\begin{equation}
m_G^{(-1)}\left(n^2\right)
=\log\left[\sec\left(\Omega_n\left(t\right)\right)+\tan\left(\Omega_n\left(t\right)\right)\right],
\end{equation}
where
\begin{equation}
t=-\frac{\pi}{2}+2\arctan\left(e^{m_G^{(-1)}(1)}\right),
\end{equation}
for every $n\in \bf R^{*}_{+}\rm$, ($t$ is a constant). Also for $k=1/2$
\begin{equation}
m_G^{(-1)}\left(4\right)=E\left[\frac{1}{2} \left(-1-\sqrt{5}+\sqrt{2 \left(5+\sqrt{5}\right)}\right),\frac{1}{2}\right].
\end{equation}
\\ 
\textbf{4)} If $n\in \bf N\rm$ and $y(x)=BR(x)$ (the Bring Radical function), then 
$$
m^{(-1)}_G(n^2)=\left(\Omega_n\circ F_1\circ b_1\right)^5+\Omega_n\circ F_1\circ b_1
$$ 
$$
b_1=3\sqrt[6]{2}\cdot{}_2F_{1}\left[\frac{1}{6},\frac{1}{3};\frac{7}{6};\frac{1}{2}\right]
$$
and $F_1$ defined from (26).

\section{Further transformations}

The function $k_i(A)$ is the inverse function of the elliptic singular modulus $k_A=k(A)$. We have  
\begin{equation}
k_i(x)=\left(\frac{K\left(\sqrt{1-x^2}\right)}{K(x)}\right)^2\textrm{, }0<x<1.
\end{equation}  
Also we define
\begin{equation}
Q_G(x):=m_{G}^{(-1)}\left(k_i(x)\right),
\end{equation} 
then
\begin{equation}
m^{(-1)}_G(r)=Q_G(k_r).
\end{equation}
From (11) we have 
$$
y\left(Q_G(k_r)\right)=R\left(e^{-\pi\sqrt{r}}\right).
$$
We are interested to find an expresion for $G$. From (10) we get
$$
5\int^{y\left(Q_G(k_r)\right)}_{0}\frac{G(t)}{t\sqrt[6]{t^{-5}-11-t^5}}dt=Q_G(k_r),
$$
or equivalently
\begin{equation}
5\int^{R(q)}_{0}\frac{G(t)dt}{t\sqrt[6]{t^{-5}-11-t^5}}dt=Q_G(k_r)
\end{equation}
and differentiating the last relation we get
$$
5\frac{G\left(R(q)\right)R'(q)}{R(q)\sqrt[6]{R(q)^{-5}-11-R(q)^{5}}}=Q'_G(k_r)\frac{dk}{dq}.
$$
Or, using (3),(5),(8) we arive to
\begin{equation}
Q'_G(k_r)=\frac{2^{1/3}}{\left(k_rk^{*}_r\right)^{2/3}}G\left(R(q)\right).
\end{equation}
But from the fact that Rogers-Ramanujan's continued fraction is algebraic function of elliptic singular moduli we have $R(q)=F(k_r)$. Hence
\begin{equation}
Q'_G(A)=\frac{2^{1/3}}{\left(A\sqrt{1-A^2}\right)^{2/3}}G\left(F(A)\right).
\end{equation}
Inverting $F$ we get\\
\begin{equation}
G(A)=2^{-1/3}\left(F_i(A)\sqrt{1-F_i(A)^2}\right)^{2/3}Q'_G\left(F_i(A)\right).
\end{equation}

The study of $F(A)$ has to reveal some interesting properties of the general function
\begin{equation}
y_s(A)=R\left(e^{-\pi\sqrt{k_i\left(s(A)\right)}}\right),
\end{equation}
where $s(A)$ is ''arbitrary'' function. Allong with $s(A)$, we attach the function $\sigma(A)$, which satisfies the condition $s'(A)=\sigma(s(A))$ and 
\begin{equation}
m_G(A)=k_i\left(s(A)\right).
\end{equation}
It holds from the definition of $y_s$ and Theorem 5 and (22): 
$$
\frac{1}{G_s\left(y_s(A)\right)}=\frac{1}{G\left(y(A)\right)}=h'_i(A)=\frac{d}{dA}\left(3\sqrt[3]{2s(A)}\cdot{}_2F_1\left[\frac{1}{3},\frac{1}{6};\frac{7}{6};s(A)^2\right]\right)=
$$
\begin{equation}
=\frac{2^{1/3}}{s(A)^{2/3}\left(1-s(A)^2\right)^{1/3}}s'(A).
\end{equation}
Hence
\begin{equation}
\frac{1}{G\left(R\left(e^{-\pi\sqrt{k_i(A)}}\right)\right)}=\frac{2^{1/3}}{A^{2/3}\left(1-A^2\right)^{1/3}}\sigma\left(A\right).
\end{equation}
Inverting the $F(x)$ function (see and Apendix) we get\\
\\
\textbf{Theorem 9.}
\begin{equation}
Q(A)=s^{(-1)}(A)
\end{equation}
and
\begin{equation}
G(A)=G_s(A)=\frac{2^{-1/3}\left(F_i(A)\sqrt{1-F_i(A)^2}\right)^{2/3}}{\sigma\left(F_i(A)\right)}.
\end{equation}
\\
\textbf{Theorem 10.}\\
Suppose $G_0(x)$ is that of (155),(65) below and $F(x)=R\left(e^{-\pi\sqrt{k_i(x)}}\right)$. Both functions are ''constant'' and algebraic. If 
\begin{equation}
5\int^{y_s(A)}_{0}\frac{G_0(t)}{\sigma\left(F_i(t)\right)}\frac{dt}{t\sqrt[6]{t^{-5}-11-t^5}}=A,
\end{equation} 
then $y_s(A)=y(A)$ is that of (58), with $s'(A)=\sigma(s(A))$ and
$$
G(t)=\frac{G_0(t)}{\sigma\left(F_i(t)\right)}.\eqno{(64.0)}
$$
\\  
\textbf{Proof.}\\
Inverting $s(A)$ we get
$$
5\int^{R\left(e^{-\pi\sqrt{k_i(A)}}\right)}_{0}\frac{G_0(t)}{\sigma\left(F_i(t)\right)t\sqrt[6]{t^{-5}-11-t^5}}dt=s_i(A).
$$
Then differentiating (the $h$ is refering to $R\left(e^{-\pi\sqrt{k_i(A)}}\right)$ function)
$$
\frac{G_0\left(R\left(e^{-\pi\sqrt{k_i(A)}}\right)\right)}{\sigma(A)}h_i'(A)=s_i'(A).
$$
Using now Corollary 2 and (60),(61), we get
$$
\frac{G_0\left(R\left(e^{-\pi\sqrt{k_i(A)}}\right)\right)}{\sigma(A)}\frac{2^{1/3}}{A^{2/3}(1-A^2)^{1/3}}=s_i'(A).
$$
Equivalently
$$
\frac{2^{-1/3}A^{2/3}(1-A^2)^{1/3}}{\sigma(A)}\frac{2^{1/3}}{A^{2/3}(1-A^2)^{1/3}}=s_i'(A)
$$
and finaly
$$
\sigma(A)s_i'(A)=1, \eqno{(64.1)}
$$
which is true hence we get (64).\\
\\
\textbf{Note.}\\
For the function $G_0(x)$ holds
\begin{equation}
G_0\left(F(x)\right)=\frac{\left(x\sqrt{1-x^2}\right)^{2/3}}{\sqrt[3]{2}}
\end{equation}
and hence $G_0$ is algebraic (see also Appendix for $G_0$).\\
\\
\textbf{Theorem 11.}\\
We set $c(A)=3\sqrt[3]{2A}\cdot{}_2F_{1}\left[\frac{1}{6},\frac{1}{3};\frac{7}{6};A^2\right]$, then
\begin{equation}
Q_G^{(-1)}\left(h(A)\right)=c_i(A).
\end{equation}
Also
\begin{equation}
\phi(A)=\frac{\sigma\left(k(A)\right)}{k'(A)},
\end{equation}
where $\phi(x)$ is that of (14). Also
$$
m^{(-1)}_G(A)=\int^{k(A)}_{0}\frac{dt}{\sigma(t)}.\eqno{(67.1)}
$$
\\
\textbf{Proof.}\\
From Theorem 5 and (35) we have $h_i\left(m_G^{(-1)}(A)\right)=b_A$ or $h_i\left(Q(A)\right)=c(A)$. Inverting we get the first result. For the second result we have
$$
m'_G(A)=\phi(m_G(A)),
$$
or from (59)
$$
k'_i(s(A))s'(A)=\phi(m_G(A)),
$$
or
$$
k'_i(s(A))s'(A)=\phi(k_i(s(A))),
$$
or
$$
\frac{k'_i(A)}{s'_i(A)}=\phi(k_i(A)),
$$
or
$$
s'_i(A)=\frac{k'_i(A)}{\phi(k_i(A))},
$$
or
\begin{equation}
s_i(A)=\int^{k_i(A)}_{0}\frac{dt}{\phi(t)}=\int^{A}_{0}\frac{dt}{\sigma(t)}.
\end{equation}
Differentiating the last relation and inverting $k_A=k(A)$ we get the result.\\
\\ 
\textbf{Note.}\\
One can see imediately that
\begin{equation}
h\left(3\sqrt[3]{2A}\cdot{}_2F_1\left[\frac{1}{3},\frac{1}{6};\frac{7}{6};A^2\right]\right)=s_i(A),
\end{equation}
which means that $h$ and $s_i$ ''generalized'' functions are esentialy the same. We are going to describe this kind of relation between generalized functions.\\
\\
\textbf{Definition 3.}\\
We say that a function $f$ is generalized, if it is not ''constant'' function.\\
\\
\textbf{Definition 4.}\\
We say that two invertible generalized functions $f,g$ are equivalent $f\equiv g$, if exist constant functions $\alpha_1(x),\beta_1(x),\gamma_1(x)$ such that
\begin{equation}
f(x)=\frac{\alpha_1\left(g\left(\beta_1(x)\right)\right)}{\gamma_1(x)}.
\end{equation}
\\
\textbf{Proposition 2.}\\
The notation $\equiv$ is an equivalence relation i.e. it has the following properties\\ 
i. \textbf{Reflection}: $f\equiv f$\\
ii. \textbf{Symmetry}: If $f\equiv g$ then $g\equiv f$\\
iii. \textbf{Transition}: If $f\equiv g$ and $g\equiv h$, then $f\equiv h$.\\
\\
\textbf{Proposition 3.}\\
We have the following equivalences of functions
\begin{equation} 
y\equiv m_G\equiv h_i\equiv Q_G^{(-1)} \equiv s
\end{equation}
\begin{equation}
y_i\equiv m_G^{(-1)}\equiv h\equiv Q_G\equiv s_i
\end{equation}
\begin{equation}
y'_i\equiv G\equiv h'\equiv\phi\equiv m_G^{(-1)}{'}\equiv Q'_G\equiv \sigma\equiv s'_i  
\end{equation}
\\
\textbf{Theorem 12.}
\begin{equation}
y'(A)=5^{-1}\sqrt[3]{2}\left(s(A)\sqrt{1-s(A)^2}\right)^{-2/3}s'(A)y(A)\sqrt[6]{y(A)^{-5}-11-y(A)^5}.
\end{equation}
\\

Set now
\begin{equation}
U(x):=256\frac{\left(1-x^2+x^4\right)^3}{x^4\left(1-x^2\right)^2}
\end{equation}
and
\begin{equation}
U_2(x)=\frac{1-\sqrt{1-U(x)^2}}{1+\sqrt{1-U(x)^2}}=16\frac{(1+14 x^2+x^4)^3}{x^2(1-x^2)^4},
\end{equation}
then we have the next\\
\\ 
\textbf{Theorem 13.}\\For an arbitrary $G(x)$ the value of $\sigma(x)$ is given from 
\begin{equation}
\sigma(x)=\frac{\left(x\sqrt{1-x^2}\right)^{2/3}}{\sqrt[3]{2}\cdot G\left(F(x)\right)}
\end{equation}
and the value of $s(x)$ from
\begin{equation}
\int^{s(x)}_{0}\frac{dt}{\sigma(t)}=x.
\end{equation}
Then the value of $y(x)$ at $x=A$ can evaluated from 
\begin{equation}
U_2\left(s\left(A\right)\right)^{1/3}Y^{5/3}=Y^2+250Y+3125,
\end{equation}    
where 
\begin{equation}
Y=y\left(A\right)^{-5}-11-y\left(A\right)^5.
\end{equation}
\\
\textbf{Notes.}\\
\textbf{I.} Assume that $G(A)=G^{*}\left(F_i(A)\right)$, where $G^{*}(A)$ known. Then from (77) we have 
\begin{equation}
\sigma(A)=2^{-1/3}\frac{\left(A\sqrt{1-A^2}\right)^{2/3}}{G^{*}(A)}.
\end{equation}
Hence
\begin{equation}
\int^{A}_{0}\frac{dt}{\sigma(t)}=\sqrt[3]{2}\int^{A}_{0}\frac{G^{*}(t)dt}{(t\sqrt{1-t^2})^{2/3}}=s_i(A).
\end{equation}
\\

Let $G(t)=\sqrt[3]{1-F_i(t)^2}$, then $s(A)=\frac{A^3}{54}$.\\
\\
\textbf{i.} If $A_0=3\sqrt[3]{6-4\sqrt{2}}$ we have 
$$
j^{1/3}=U_2\left(s\left(A_0\right)\right)^{1/3}=12.
$$
Then equation (79) becomes
$$
12Y^{5/3}=Y^2+250Y+3125
$$
and
$$
Y=125(2+\sqrt{5}).
$$
Hence
$$
y(A_0)=y\left(\sqrt[3]{6-4\sqrt{2}}\right)=\frac{2}{\sqrt{2 \left(5+\sqrt{5}\right)}+\sqrt{5}+1}.
$$
Finaly the function $y(A)$ which is a solution of
$$
5\int^{y(A)}_{0}\frac{\sqrt[3]{1-F_i(t)^2}}{t\sqrt[6]{t^{-5}-11-t^5}}dt=A,
$$
is
$$
y(A)=R\left(e^{-\pi\sqrt{k_i\left(A^3/54\right)}}\right)
$$
and can be determined in closed form up to a 6th degree poynomial equation (that of (79)).\\
\textbf{ii.} If 
$$
A_0=3 \sqrt[3]{66+48 \sqrt{2}-8 \sqrt{140+99 \sqrt{2}}},
$$
then
$$
j^{1/3}=U_2\left(s(A_0)\right)^{1/3}=66
$$
and equation (79) becomes
$$
66Y^{5/3}=Y^2+250 Y+3125,
$$
with solution
$$
Y=\frac{125}{2} \left(1147+513\sqrt{5}+\sqrt{2630810+1176534 \sqrt{5}}\right).
$$
Hence
$$
y\left(3 \sqrt[3]{66+48 \sqrt{2}-8 \sqrt{140+99 \sqrt{2}}}\right)=\sqrt[5]{\frac{-11-Y+\sqrt{125+22Y+Y^2}}{2}}.
$$
\\
\textbf{II.} If for a function $y(A)$ we know $G$, then from (77) we have 
$$
\int^{A}_{0}\frac{dt}{\sigma(t)}=s_i(A).\eqno{:(eq1)}
$$
Knowing $s(A)$ we solve
$$
\sqrt[3]{16}\frac{1+14 s(A)^2+s(A)^4}{s(A)^{2/3} \left(1-s(A)^2\right)^{4/3}}\cdot Y^{5/3}=Y^2+250Y+3125\eqno{:(eq2)} 
$$
and we get that 
$$
y(A)=y_s(A)=\sqrt[5]{\frac{-11-Y+\sqrt{125+22Y+Y^2}}{2}}.\eqno{:(eq3)}
$$ 
Hence we find the closed form of $y(A)$ in (64) (and hence to the problem (10)) up to the inverting of the integral of $(eq1)$ and solving the sextic equation $(eq2)$.\\
\\
\textbf{i)} Suppose that $\sigma(A)=A+1$, then $s(A)=e^A-1$. This case coresponds to
$$
5\int^{y(x)}_{0}\frac{G_0(t)}{F_i(t)+1}\frac{dt}{t\sqrt[6]{t^{-5}-11-t^5}}=x
$$ 
and is solvable up to the sextic equation $(eq2)$. Also
$$
y(A)=R\left(e^{-\pi\sqrt{k_i(e^A-1)}}\right).
$$
\\
\textbf{ii)} Another example is with $\sigma(A)=1$. This leads to $s(A)=A$ and coresponds to
$$
5\int^{y(x)}_{0}\frac{G_0(t)}{t\sqrt[6]{t^{-5}-11-t^5}}dt=x.
$$
The first derivative according to Theorem 12 is
$$
y'(A)=5^{-1}\sqrt[3]{2}\left(A\sqrt{1-A^2}\right)^{-2/3}y(A)\sqrt[6]{y(A)^{-5}-11-y(A)^5}
$$
and
$$
y(A)=R\left(e^{-\pi\sqrt{k_i(A)}}\right).
$$
For more details see section Applications.\\
\\
\textbf{iii)} For $\sigma(A)=1/A$ we get $s(A)=\sqrt{2A}$, hence
$$
5\int^{y(x)}_{0}\frac{G_0(t)F_i(t)dt}{t\sqrt[6]{t^{-5}-11-t^5}}=x
$$
and
$$
y'(A)=\frac{5^{-1}}{\sqrt{2}\left(1-2A\right)^{1/3}A^{5/6}}y(A)\sqrt[6]{y(A)^{-5}-11-y(A)^5},
$$
with
$$
y(A)=R\left(e^{-\pi\sqrt{k_i\left(\sqrt{2A}\right)}}\right).
$$
\\ 
\textbf{iv)} If $\sigma(A)=\sqrt{1-A^2}\sqrt{1-k A^2}$, then $s(A)=\textrm{sn}(A,k)$ and the solution $y(x)$ of
$$
5\int^{y(x)}_{0}\frac{G_0(t)}{\sqrt{1-F_i(t)^2}\sqrt{1-k F_i(t)^2}}\frac{dt}{t\sqrt[6]{t^{-5}-11-t^5}}=x,
$$
is given from $(eq2)$ and $(eq3)$. Esentialy the function $y(x)=y_k(x)=y(x,k)$ is algebraic function of $s(A)=\textrm{sn}(A,k)$ and hence double periodic elliptic function. Also
$$
\frac{y'_k(A)^2}{1-kF_i\left(y_k(A)\right)^2}=\frac{y'_l(A)^2}{1-lF_i\left(y_l(A)\right)^2}=C(A)
$$
and
$$
y'_k(A)=5^{-1}\sqrt[3]{2}\frac{\textrm{cn}(A,k)^{1/3}\textrm{dn}(A,k)}{\textrm{sn}(A,k)^{2/3}}y_k(A)\sqrt[6]{y_k(A)^{-5}-11-y_k(A)^5},
$$
with
$$
y(A)=y_k(A)=R\left(e^{-\pi\sqrt{k_i\left(sn(A,k)\right)}}\right).
$$
\textbf{III.} Another notation but not so detailed can found using (58),(77),(78) and relation $R(q)=F(k_r)$. We have
\begin{equation}
y_i(A)=\sqrt[3]{2}\int^{F_i(A)}_{0}\frac{G(F(t))}{\left(t\sqrt{1-t^2}\right)^{2/3}}dt.
\end{equation}
As application we set $B(x,a,b)=\int^{x}_{0}t^{a-1}(1-t)^{b-1}dt$ to be the incoplete beta function. Then if
$$
G\left(F(A)\right)=\frac{1}{\sqrt[3]{2}}\left(A\sqrt{1-A^2}\right)^{2/3}\left(A-A^2\right)^{a-1},
$$
we have
\begin{equation}
y\left(B(x,a,a)\right)=F(x)
\end{equation}
and (see [13]) the solution of
\begin{equation}
\frac{B\left(1-\beta_r,a,a\right)}{B\left(\beta_r,a,a\right)}=r\textrm{, }r>0\textrm{, }0<\beta_r<1,
\end{equation}
is equivalent to
\begin{equation}
B\left(\beta_r ,a,a\right)=\frac{\Gamma(a)^2}{\Gamma(2a)(r+1)}.
\end{equation}
Hence we get
\begin{equation}
y\left(\frac{\Gamma(a)^2}{\Gamma(2a)(r+1)}\right)=F\left(\beta_r\right).
\end{equation}
\\
\textbf{IV.}
Also from
$$
\int^{A}_{0}\frac{dt}{\sigma(t)}=s_i(A),
$$
we have
$$
\int^{k_r}_{0}\frac{dt}{\sigma(t)}=s_i\left(k_r\right).
$$
Hence from (59): $s_i(A)=m_G^{(-1)}(k_i(A))$ or equivalently $s_i(k_r)=m_G^{(-1)}(r)$ and we get
\begin{equation}
\int^{k_r}_0\frac{dt}{\sigma(t)}=m_G^{(-1)}(r).
\end{equation}
\textbf{V.} If $G(F(A))$ is polynomial
\begin{equation}
G\left(F(x)\right)=\sum^{M}_{n=0}c_nx^n,
\end{equation} 
then using
\begin{equation}
\int^{A}_{0}\frac{t^n}{(t\sqrt{1-t^2})^{2/3}}dt=3\frac{A^{n+1/3}}{3n+1}{}_2F_{1}\left[\frac{1}{3},\frac{3n+1}{6};\frac{3n+7}{6};A^2\right],
\end{equation}
we get
\begin{equation}
y\left(3\sqrt[3]{2}\sum^{M}_{n=0}c_n\frac{A^{n+1/3}}{3n+1}\cdot {}_2F_{1}\left[\frac{1}{3},\frac{3n+1}{6};\frac{3n+7}{6};A^2\right]\right)=F(A),
\end{equation} 
or if someone preferes
\begin{equation}
3\sqrt[3]{2}\sum^{M}_{n=0}c_n\frac{A^{n+1/3}}{3n+1}\cdot {}_2F_{1}\left[\frac{1}{3},\frac{3n+1}{6};\frac{3n+7}{6};A^2\right]+c=y_i\left(F(A)\right).
\end{equation}
and concequently:\\
\\
\textbf{Theorem 14.}\\
If
\begin{equation}
\sigma(A)=\frac{(A\sqrt{1-A^2})^{2/3}}{\sqrt[3]{2}\sum^{M}_{m=0}c_mA^m},
\end{equation}
then
\begin{equation}
y_i(A)=3\sqrt[3]{2}\sum^{M}_{n=0}c_n\frac{F_i(A)^{n+1/3}}{3n+1}\cdot {}_2F_{1}\left[\frac{1}{3},\frac{3n+1}{6};\frac{3n+7}{6};F_i(A)^2\right]+c.
\end{equation}
\\

Assume function $G$ as in (89), then from the fact that $F$ is algebraic, there exists coefficients $a_{n,l}$ such that
\begin{equation}
\sum^{N}_{n,l=0}a_{kl}G(x)^nx^l=0.
\end{equation}  
But then also
$$
\sum^{N}_{n,l=0}a_{nl}G\left(F(x)\right)^nF(x)^l=0
$$
and
\begin{equation}
\sum^{N}_{n,l=0}a_{nl}\left(\sum^{M}_{m=0}c_mx^m\right)^nF(x)^l=0.
\end{equation}
Finaly
\begin{equation}
\sum^{N}_{n,l=0}a_{nl}\left(\sum^{M}_{m=0}c_m(k_r)^m\right)^nR(q)^l=0.
\end{equation}
Hence knowing $c_m$ we find $a_{nl}$ from relation (97) and equation 
$$
x^2 \left(1-x^2\right)^4 \left(y^{20}-228 y^{15}+494 y^{10}+228 y^5+1\right)^3+
$$
\begin{equation}
+16 \left(x^4+14 x^2+1\right)^3 y^5 \left(y^{10}+11 y^5-1\right)^5=0,
\end{equation}
(which is Klein formula for the icosahedron $y=F(x)$). That is the $a_{nl}$ of (95) can be found from that of $c_m$, by equating coefficients of the identity: 
$$
\sum^{N}_{n,l=0}a_{nl}\left(\sum^{M}_{m=0}c_mx^m\right)^ny^l=x^2 \left(1-x^2\right)^4 \left(y^{20}-228 y^{15}+494 y^{10}+228 y^5+1\right)^3+
$$
\begin{equation}
+16 \left(x^4+14 x^2+1\right)^3 y^5 \left(y^{10}+11 y^5-1\right)^5.
\end{equation}
If hapens $G(F(x))=\psi(x)$ be more complicated, for example algebraic, then we solve the equation $G=\psi(x)$ with respect to $x$, $x=\psi^{(-1)}(G)=\psi_i(G)$ and
\begin{equation}
\sum^{60}_{n,l=0}A_{nl}\left(\psi_i(G(x))\right)^nx^l=0,
\end{equation}
where the $A_{nl}$ are that of Klein's equation 
$$
\sum^{60}_{n,l=0}A_{nl}x^ny^l=x^2 \left(1-x^2\right)^4 \left(y^{20}-228 y^{15}+494 y^{10}+228 y^5+1\right)^3+
$$
\begin{equation}
+16 \left(x^4+14 x^2+1\right)^3 y^5 \left(y^{10}+11 y^5-1\right)^5.
\end{equation}
Hence in general:\\ 
\\
\textbf{Theorem 15.}\\
If $G(F(x))=\psi(x)$ is known resonable function (polynomial algebraic etc...), then the minimal equation for $G$ is (100), with coefficients $A_{nl}$ that of (101).\\
\\
\textbf{Example.}\\
For $G(F(x))=2^{-1/3}(x^2-x^4)^{1/3}$, we get $Q'_G(x)=1$, hence $s_i(x)=x$ and $y(x)=F(x)$, where $G(x)=G_0(x)$ is solution of 
\begin{equation}
\sum^{60}_{n,l=0}A_{nl}\left(\sqrt{\frac{1-\sqrt{1-8G(x)^3}}{2}}\right)^nx^l=0.
\end{equation}
The coefficients $A_{nl}$ are that of (101).

\section{Solution of General Equations and Inversion}

Consider a polynomial $w(x)$ and the equation
\begin{equation}
w(x)=\lambda
\end{equation}
Set $G(x)=5^{-1}xw'(x)\sqrt[6]{x^{-5}-11-x^5}$, then  
\begin{equation}
m_G^{(-1)}(r)=w\left(R\left(e^{-\pi\sqrt{r}}\right)\right)
\end{equation}
and the solution of (103) is $x=R\left(e^{-\pi\sqrt{r}}\right)$, where $r=m_{G}(\lambda)$.\\

Taking the derivatives with respect to $A$ in $w(f^{(-1)}(A))=A$, we lead to
\begin{equation}
w{'}\left(f^{(-1)}\left(A\right)\right)f^{(-1)}{'}(A)=1,
\end{equation}
or
\begin{equation}
f^{(-1)}{'}\left(A\right)=\frac{1}{w{'}\left(f^{(-1)}\left(A\right)\right)}.
\end{equation}
Hence if we set $W(A)=\frac{1}{w{'}\left(A\right)}$,
then clearly 
\begin{equation}
f^{(-1)}{'}\left(A\right)=W\left(f^{(-1)}\left(A\right)\right).
\end{equation}
Setting where $f^{(-1)}(A)=y(A)$, with $y(A)$ that of (9) and (10) and taking the derivatives with respect to $A$ we have 
$$
5\frac{G\left(f^{(-1)}(A)\right)f^{(-1)}{'}(A)}{f^{(-1)}(A)\sqrt[6]{\left(f^{(-1)}(A)\right)^{-5}-11-\left(f^{(-1)}(A)\right)^{5}}}=
$$
\begin{equation}
=5\frac{G\left(f^{(-1)}(A)\right)W\left(f^{(-1)}(A)\right)}{f^{(-1)}(A)\sqrt[6]{\left(f^{(-1)}(A)\right)^{-5}-11-\left(f^{(-1)}(A)\right)^{5}}}=1. 
\end{equation}
After inverting $f^{(-1)}(x)$ we have
$$
G(x)W(x)=5^{-1}x\sqrt[6]{x^{-5}-11-x^5},
$$
or equivalently
\begin{equation}
G(x)=5^{-1}xw'(x)\sqrt[6]{x^{-5}-11-x^5}
\end{equation}
and it will be
\begin{equation}
f^{(-1)}(A)=R\left(e^{-\pi\sqrt{m_G(A)}}\right).
\end{equation}
From the above we can state the following theorem\\ 
\\
\textbf{Theorem 16.}\\
The equation $w(y(x))=x$ have solution 
\begin{equation}
y(x)=R\left(e^{-\pi\sqrt{m_G(x)}}\right),
\end{equation}
with $G$ that of relation (109) and $m_G(A)$ as defined in (9).\\
\\
\textbf{Example.}\\
Let $\rho_1=\frac{1}{2} \left(11-5 \sqrt{5}\right)$, $\rho_2=\frac{1}{2} \left(11+5 \sqrt{5}\right)$ and
consider the equation
\begin{equation}
6\frac{x^{a+\frac{11}{6}}}{6 a+11} F_{Ap}\left(\frac{6 a+11}{30},\frac{1}{6},\frac{1}{6},\frac{6 a+41}{30},\rho_1 x^5,\rho_2 x^5\right)=A,
\end{equation}
where $a$ is parameter. The function $G$ is $G(x)=\frac{x^{a+1}}{5}$, and 
\begin{equation}
x=R\left(e^{-\pi\sqrt{m_G(A)}}\right),
\end{equation}
where the $m_G(x)$ is given from
\begin{equation}
x=\frac{\pi}{5}\int^{\infty}_{\sqrt{m_G(x)}}\eta\left(it/2\right)^4R\left(e^{-\pi t}\right)^{a+1}dt.
\end{equation}
Hence if
$$
g_a(x)=6 \frac{x^{a+\frac{11}{6}}}{6 a+11} F_{Ap}\left(\frac{6 a+11}{30},\frac{1}{6},\frac{1}{6},\frac{6 a+41}{30},\rho_1 x^5,\rho_2 x^5\right),
$$
then
\begin{equation}
g_a\left(R\left(e^{-\pi A}\right)\right)=\frac{\pi}{5}\int^{\infty}_{A}\eta\left(it/2\right)^4R\left(e^{-\pi t}\right)^{a+1}dt.
\end{equation}
\\
\textbf{Theorem 17.}\\
Given the equation $P(x)=a:(\epsilon)$, the inverse of $P(x)$ is $y(x)$, then $(\epsilon)$ is equivalent to
\begin{equation}
\int^{F_i(x)}_{0}\frac{dt}{\sigma(t)}=a.
\end{equation}
\\
\textbf{Proof.}\\
Easy\\
\\
\textbf{Example.}\\
If $\sigma(x)=x+1$, then 
$$
P(x)=5\int^{x}_{0}\frac{G_0(t)}{F_i(t)+1}\frac{dt}{t\sqrt[6]{t^{-5}-11-t^5}}
$$
and the equation $P(x)=a$ have solution $x$ such that
$$
\int^{F_i(x)}_{0}\frac{dt}{t+1}=a,
$$ 
or equivalently $x=F\left(e^a-1\right)$.\\
\\

In general we have the next formula\\
\\
If $|x|<1$ then
\begin{equation}
\int^{x}_{0}\frac{f(-q)^5R(q)g'\left(R(q)\right)}{f(-q^5)}dq=5g\left(R\left(x\right)\right),
\end{equation}
which is consequence of the next identity 
\begin{equation}
\frac{R'(q)}{R(q)}=\frac{f(-q)^5}{5qf(-q^5)}.
\end{equation}
Relation (118) was given by Ramanujan (see [3]). We know that
$$
\int^{q_1}_{q_2}f(-q)^4q^{-5/6}R(q)^{5\nu}dq=-\pi\int^{\sqrt{r_1}}_{\sqrt{r_2}}\eta(it/2)^4R(e^{-\pi t})^{5\nu}dt
$$
and (see [13]):
$$
C(\nu):=\int^{1}_{0}f(-q)^4q^{-5/6}R(q)^{5\nu}dq=
$$
\begin{equation}
=\Gamma\left(\frac{5}{6}\right)\left(\frac{11+5 \sqrt{5}}{2}\right)^{-\frac{1}{6}-\nu} \frac{\Gamma\left(\frac{1}{6}+\nu\right)}{\Gamma(1+\nu)} {}_2F_1\left(\frac{1}{6},\frac{1}{6}+\nu;1+\nu;\frac{11-5 \sqrt{5}}{11+5\sqrt{5}}\right),
\end{equation}
where $\nu\geq0$. Hence\\ 
\\
\textbf{Theorem 18.}\\
If $G(x)$ is a polynomial (or analytic function when $n\rightarrow+\infty$) of the form 
\begin{equation}
G(x)=\sum^{n}_{m=0}a_nx^{p_{m}},
\end{equation} 
with $p_m-$positive reals, ($\lim p_m=+\infty$) and $R\left(1\right)=\frac{\sqrt{5}-1}{2}$, then\\
\textbf{i)}
\begin{equation}
m_G^{(-1)}(0)=\pi\int^{\infty}_{0}\eta(it/2)^4G\left(R(e^{-\pi t})\right)dt=\sum^{n}_{m=0}a_m C\left(5^{-1}p_m\right).
\end{equation}
\textbf{ii)} If $y$ is a smooth function and $G$ is of the form (120) then the equation 
\begin{equation}
y(x)=\frac{\sqrt{5}-1}{2},
\end{equation}
have a solution
\begin{equation}
x=x_0=\sum^{n}_{m=0}a_mC(5^{-1}p_m).
\end{equation} 
\\
\textbf{Example.}\\
Let $G(x)=e^{-x}-1$. Then $a_m=\frac{(-1)^m}{m!}$, $m=1,2,\ldots$ and
\begin{equation}
5\int^{y(x)}_{0}\frac{e^{-t}-1}{t\sqrt[6]{t^{-5}-11-t^5}}dt=x
\end{equation}
and the solution of $y(x)=\frac{\sqrt{5}-1}{2}$ is
\begin{equation}
x=\sum^{\infty}_{m=1}\frac{(-1)^m}{m!}C\left(\frac{m}{5}\right).
\end{equation}
Note here that if $G$ has finite expansion (120) the result (123) becomes more meaningful since the evaluation of the root of (122) by hypergeometric functions is beter than the integral:
\begin{equation}
5\int^{\frac{\sqrt{5}-1}{2}}_{0}\frac{G(t)}{t\sqrt[6]{t^{-5}-11-t^5}}dt.
\end{equation}    
\\

Assume now 
$$
G^{*}(t):=G(t\alpha)\sqrt[6]{\frac{t^{-5}-11-t^5}{(\alpha t)^{-5}-11-(\alpha t)^5}}.\eqno{(126.1)}
$$
If $y^{*}(A)$ coresponds to $G^{*}(t)$, one can easily see that
$$
y^{*}(A)=\frac{y(A)}{\alpha}.\eqno{(126.2)}
$$
Hence we have the next theorem which is generalization of Theorem 18:\\
\\
\textbf{Theorem 18.1}\\
Assume the function $G(t)$ is given near the origin by (under certain converging conditions):
$$
G(t)=\sum^{\infty}_{m=0}a_mt^{p_m},\eqno{(126.3)}
$$
where $p_m$ is any increasing sequence of positive real numbers with $\lim p_m=+\infty$. Then the function $y(A)$ defined as 
$$
5\int^{y(A)}_{0}\frac{G(t)}{t\sqrt[6]{t^{-5}-11-t^5}}dt=A,\eqno{(126.4)}
$$
have the following property: Every equation of the form
$$
y(x)=\alpha \frac{\sqrt{5}-1}{2},\eqno{(126.5)}
$$
have solution $x$ such that
$$
x=\sum^{\infty}_{m=0}a_m^{*}(\alpha)C\left(5^{-1}p^{*}_m\right),\eqno{(126.6)}
$$
where $a^{*}_m(\alpha)$ is such that
$$
G(x\alpha)\sqrt[6]{\frac{x^{-5}-11-x^5}{(\alpha x)^{-5}-11-(\alpha x)^5}}=\sum^{\infty}_{m=0}a^{*}_{m}(\alpha)x^{p^{*}_m}.\eqno{(126.7)}
$$
Hence if we set $\xi^{-1}=\frac{\sqrt{5}-1}{2}$, then the inverse of $y(A)$ is
$$
y_i(A)=\sum^{\infty}_{m=0}a^{*}_m(\xi A)C\left(5^{-1}p^{*}_m\right),\eqno{(126.8)}
$$
provited the convergence of (126.1),(126.3),(126.4),(126.6),(126.7),(126.8).\\
\\

Set 
$$
G=G_1(t)=\frac{t\sqrt[6]{t^{-5}-11-t^5}}{5\sqrt{1-t^2}\sqrt{1-k^2t^2}}.
$$
Then
\begin{equation}
5\int^{y}_{0}\frac{G_1(t)}{t\sqrt[6]{t^{-5}-11-t^5}}dt=\int^{y}_{0}\frac{1}{\sqrt{1-t^2}\sqrt{1-k^2 t^2}}dt=A.
\end{equation} 
Hence (see [4]):
$$
y=\textrm{sn}(A,k)=R\left(e^{-\pi\sqrt{m_G}}\right)\textrm{ and }A=h\left(b_{m_G}\right),
$$
where
$$
b_{m_G}=3\sqrt[3]{2k_{m_G}}\cdot{}_2F_{1}\left[\frac{1}{6},\frac{1}{3};\frac{7}{6};k_{m_G}^2\right]
$$
and
$$
h^{(-1)}(x)=\int^{x}_{0}\frac{du}{G_1\left(\textrm{sn}\left(u,k_r\right)\right)}=5\int^{x}_{0}\frac{\textrm{dn}(u)\textrm{cn}(u)}{\textrm{sn}(u)\sqrt[6]{\textrm{sn}(u)^{-5}-11-\textrm{sn}(u)^5}}du.
$$
Hence if $y(A)$ is defined by (9),(10) and $\textrm{sn}(u)=\textrm{sn}(u,k)$, $\textrm{dn}(u)=\textrm{dn}(u,k)$, then\\
\\
\textbf{Theorem 19.}
$$
5\int^{E\left[\arcsin(y(A)),k\right]}_{0}\frac{\textrm{dn}(u)\textrm{cn}(u)}{\textrm{sn}(u)\sqrt[6]{\textrm{sn}(u)^{-5}-11-\textrm{sn}(u)^5}}du=
$$
\begin{equation}
=3\sqrt[3]{2k_{m_G}}\cdot{}_2F_{1}\left[\frac{1}{6},\frac{1}{3};\frac{7}{6};k_{m_G}^2\right],
\end{equation}
where $k$ is independent parameter, $m_G=m_G(A)$ and $E$ denotes the incomplete elliptic integral of the first kind, $y$ is the function defined in (9),(10),(11).\\
\\ 

Inverting the above integral we get a formula for the Rogers-Ramanujan continued fraction:\\
Set 
\begin{equation}
5\int^{H_o(x)}_{0}\frac{\textrm{dn}(u)\textrm{cn}(u)}{\textrm{sn}(u)\sqrt[6]{\textrm{sn}(u)^{-5}-11-\textrm{sn}(u)^5}}du=x,
\end{equation}
then
\begin{equation}
R(q)=\textrm{sn}\left(H_o(b_r),k\right).
\end{equation} 
For the function $\textrm{sn}$ we have $G(t)=G_1(t)$ and if the equation
$$
\textrm{sn}(x,k)=a,
$$
have $m_1$ such that $R\left(e^{-\pi\sqrt{m_1}}\right)=a$. 
The solution is $x=x_1$:
$$
x_1=E\left[\arcsin(a),k\right]=\pi\int^{+\infty}_{\sqrt{m_1}}\eta\left(it/2\right)^4G_1\left(R\left(e^{-\pi t}\right)\right)dt.
$$
Taking derivatives in (129) we get 
\begin{equation}
H_o{'}\left(H_o^{(-1)}(x)\right)=\frac{1}{H_o^{(-1)}{'}(x)}=\frac{\textrm{sn}(x)\sqrt[6]{\textrm{sn}(x)^{-5}-11-\textrm{sn}(x)^5}}{5\textrm{cn}(x)\textrm{dn}(x)}.
\end{equation}
Hence
\begin{equation}
H_o^{(-1)}{'}\left(E\left[\arcsin\left(x\right),k^2\right]\right)=\frac{5\sqrt{1-k^2x^2}\sqrt{1-x^2}}{x\sqrt[6]{x^{-5}-11-x^5}}.
\end{equation}
According to (129) and (130), the function $H_o=H_o(x,m)$, takes special values
$$
H_o(b_r)=H_o(b_r,m)=E(\arcsin(R(q)),m),
$$
where $q=e^{-\pi\sqrt{r}}$, $r>0$ for all $0<m<1$.\\
By this way $H_o$ can evaluated with known functions
$$
H_o(A,k)=E\left(\arcsin\left(R\left(e^{-\pi\sqrt{m(A)}}\right)\right),k\right),\eqno{(132.1)}
$$
with $m(A)$ that of (23),(24).

\section{More integrals}

Let $F_1$ be the function introduced in the above sections i.e. $F_1(x)=R\left(e^{-\pi\sqrt{m(x)}}\right)$, then $F_1$ is such that
\begin{equation}
F'_1(x)=5^{-1}F_1\left(x\right)\sqrt[6]{\left(F_1(x)\right)^{-5}-11-\left(F_1(x)\right)^5}.
\end{equation}
Also $F_1^{(-1)}$ is a specific Appell function 
\begin{equation}
F_{Ap}[a,b_1,b_2,c,x,y]:=\sum^{\infty}_{m=0}\sum^{\infty}_{n=0}\frac{(a)_{m+n}}{(c)_{m+n}}\frac{(b_1)_m (b_2)_{n}}{m! n!}x^m y^n .
\end{equation}
More precicely
\begin{equation}
F_1^{(-1)}(x)=6x^{5/6}F_{Ap}\left[\frac{1}{6},\frac{1}{6},\frac{1}{6},\frac{7}{6},\frac{-2x^5}{11+5\sqrt{5}},\frac{-2x^5}{11-5\sqrt{5}}\right].
\end{equation}
Set also $t=t(w)$ such that
\begin{equation}
t(w)=F_1\left[(-1)^{m+1}a^{m-1}D^{-m+1/2}B\left(\frac{-b+\sqrt{D}-2aw}{2\sqrt{D}},1-m,1-m\right)\right],
\end{equation}
where $B(x,a,b)=\int^{x}_{0}t^{a-1}(1-t)^{b-1}dt$, $D=b^2-4ac$.\\
Also let 
\begin{equation}
U(a,b,c;m;x):=(-1)^{m+1}a^{m-1}D^{-m+1/2}B\left(\frac{-b+\sqrt{D}-2ax}{2\sqrt{D}},1-m,1-m\right).
\end{equation}
Then 
$$
\frac{d}{dx}U(a,b,c;m;x)=\frac{1}{(ax^2+bx+c)^m}\Leftrightarrow
$$
\begin{equation} \int^{B}_{A}\frac{dt}{(at^2+bt+c)^m}=U(a,b,c;m;B)-U(a,b,c;m;A).
\end{equation}
Hence we get 
$$
5\int^{t_1}_{t_0}\frac{G(t)dt}{t\sqrt[6]{t^{-5}-11-t^5}}=
$$
\begin{equation}
=\int^{w_1}_{w_0} \frac{G\left(F_1\left[(-1)^{m+1}a^{m-1}D^{-m+1/2}B\left(\frac{-b+\sqrt{D}-2aw}{2\sqrt{D}},1-m,1-m\right)\right]\right)}{(aw^2+bw+c)^m}dw ,
\end{equation}
where $w_1,w_0,t_0,t_1$ are such that
\begin{equation}
U_i\left(a,b,c;m;F_1^{(-1)}\left(t_1\right)\right)=w_1\textrm{ and }U_i\left(a,b,c;m;F_1^{(-1)}\left(t_0\right)\right)=w_0.
\end{equation}

Suppose that we wish to evaluate the integral
\begin{equation}
I:=\int^{w_1}_{w_0}\frac{f(w)}{\sqrt{aw^2+bw+c}}dw.
\end{equation}
Then easily from (139) with $m=1/2$ 
$$
I=5\int^{F_1\left(U\left(w_1\right)\right)}_{F_1\left(U\left(w_0\right)\right)}\frac{f\left(U_i\left(F_1^{(-1)}(t)\right)\right)}{t\sqrt[6]{t^{-5}-11-t^5}}dt=
$$
\begin{equation}
=y_i\left(F_1\left(U\left(w_1\right)\right)\right)-y_i\left(F_1\left(U\left(w_0\right)\right)\right)=h\left(U(w_1)\right)-h\left(U(w_0)\right),
\end{equation}
where $y(x)=R\left(e^{-\pi\sqrt{m_G(x)}}\right)$ and
\begin{equation}
G(x)=f\left(U_i\left(F_1^{(-1)}(x)\right)\right)
\end{equation}
and 
\begin{equation}
U_i(x):=\frac{1}{2a} \left(\sqrt{b^2-4 a c}\cdot \sinh \left(\sqrt{a} x\right)-b\right).
\end{equation}
Until now we have evaluated $G(x)$. From Theorem 4 we have
$$
h'\left(5\int^{x}_{0}\frac{dt}{t\sqrt[6]{t^{-5}-11-t^5}}\right)5\frac{1}{\sqrt[6]{x^{-5}-11-x^5}}=5\frac{G(x)}{x\sqrt[6]{x^{-5}-11-x^5}}.
$$
Hence
\begin{equation}
h'(x)=G\left(F_1(x)\right)=f\left(U_i(x)\right),
\end{equation}
which is a resonable equation to find $h$. In general we can state the following\\
\\
\textbf{Theorem 20.}\\
We have
\begin{equation}
\int^{w_1}_{w_0}\frac{f(w)}{(aw^2+bw+c)^m}dw=h\left(U\left(w_1\right)\right)-h\left(U\left(w_0\right)\right),
\end{equation}
where
\begin{equation}
G(x)=f\left(U_i\left(F_1^{(-1)}(x)\right)\right).
\end{equation}
\\
\textbf{Theorem 20.1}\\
Assume that
$$
G(x)=\sum^{\infty}_{n=0}\frac{G^{(n)}(0)}{n!}x^n.\eqno{(147.1)}
$$
Then
$$
h(A)=Q_G(c_i(A))=5F_1(A)^{5/6}\times
$$
$$
\times\sum^{\infty}_{n=0}\frac{G^{(n)}(0)}{n!}\frac{F_1(A)^{n}}{5/6+n}F_{Ap}\left[\frac{1}{6}+\frac{n}{5},\frac{1}{6},\frac{1}{6},\frac{7}{6}+\frac{n}{5},\rho_1 F_1(A)^5,\rho_2 F_1(A)^5\right].\eqno{(147.2)}
$$
\\
\textbf{Proof.}\\
See section Applications paragraph 6.2.\\
\\
\textbf{Notes.}\\
More general if
$$
G(x)=\sum^{\infty}_{n=0}G_nx^{p_n},\eqno{(147.3)}
$$
where $p_n$ is increasing sequence of positive real numbers, with $\lim p_n=+\infty$, then
$$
h(A)=Q_G(c_i(A))=5F_1(A)^{5/6}\times
$$
$$
\times\sum^{\infty}_{n=0}G_n\frac{F_1(A)^{p_n}}{5/6+p_n}F_{Ap}\left[\frac{1}{6}+\frac{p_n}{5},\frac{1}{6},\frac{1}{6},\frac{7}{6}+\frac{p_n}{5},\rho_1 F_1(A)^5,\rho_2 F_1(A)^5\right].\eqno{(147.4)}
$$
\\
\textbf{Theorem 21.}\\
It is $G\left(F_1(x)\right)=f\left(U_i(x)\right)$ and
\begin{equation}
Q_G(A)=h(c_A)=\int^{c_A}_{0}f(U_i(t))dt=\int^{c_A}_{0}G(F_1(t))dt. 
\end{equation} 
\\ 
\textbf{Proof.}\\
From (145) we have $h'(b_A)=G(R(q))=f(U_i(b_A))$ inverting $k(A)$ we get $h'(c_A)=G(F(A))=f(U_i(c_A))$. Hence  $h'(c_A)c'_A=f(U_i(c_A))c'_A$ and integrating, $h(c_A)=Q_G
(A)=\int^{c_A}_{0}f(U_i(t))dt=\int^{c_A}_{0}G(F_1(t))dt$. Hence we get the result.\\  
\\
\textbf{Example.}\\
If $G(t)=\sqrt[6]{t}$, then from (167) below
$$
Q_G(A)=\int^{c_A}_{0}\sqrt[6]{F_1(t)}dt=
$$
$$
=5F_1(c_A)F_{Ap}\left[\frac{1}{5},\frac{1}{6},\frac{1}{6},\frac{6}{5},\frac{11-5\sqrt{5}}{2}F_1(c_A)^5,\frac{11+5\sqrt{5}}{2}F_1(c_A)^5\right].
$$
Hence $y(A)=R\left(e^{-\pi\sqrt{k_i\left(Q_G^{(-1)}(A)\right)}}\right)$ and holds the following semi-algebraic relation for the function $y(x)$:
$$
y\left(5F_1(c_A)F_{Ap}\left[\frac{1}{5},\frac{1}{6},\frac{1}{6},\frac{6}{5},\frac{11-5\sqrt{5}}{2}F_1(c_A)^5,\frac{11+5\sqrt{5}}{2}F_1(c_A)^5\right]\right)=F(A).
$$
\\
\textbf{Theorem 22.}\\
Given $G$ there holds the relation
\begin{equation}
y_i\left(G_i(x)\right)=xF_1^{(-1)}(G_i(x))-\int^{x}_{0}F_1^{(-1)}\left(G_i(t)\right)dt,
\end{equation}
where $F_1^{(-1)}$ is the Appell function of (135).\\
\\
\textbf{Proof.}\\
Integration by parts.\\
\\
\textbf{Example.}\\
Suppose $G(x)=\log\left(F_1^{(-1)}(x)+1\right)$, then $F_1^{(-1)}(G_i(x))=e^x-1$ hence
$$
y_i\left(F_1\left(e^x-1\right)\right)=e^x(x-1)
$$
and $h(x)=\int^{x}_{0}\log(t+1)dt=(x+1)\log(x+1)-x$.\\

Also from Theorems 17 and 19 we get\\
\\ 
\textbf{Theorem 23.}\\
Let $\xi^{-1}=\sqrt[5]{\frac{-11+5\sqrt{5}}{2}}$ and 
\begin{equation}
U_i\left(a,b,c;m;F_1^{(-1)}\left(\xi^{-1}\right)\right)=p_1\textrm{ and }U_i\left(a,b,c;m;F_1^{(-1)}(0)\right)=p_0,
\end{equation}
then
\begin{equation}
\int^{p_1}_{p_0} \frac{G\left(F_1\left[U(a,b,c;m;w)\right]\right)}{(aw^2+bw+c)^m}dw=\sum^{\infty}_{n=0}\frac{G^{(n)}(0)}{n!}C\left(\frac{n}{5}\right).
\end{equation}
\\
\textbf{Corollary 1.}
$$
5\int^{p_1}_{p_0} \frac{F_1\left[U\left(a,b,c;m;w\right)\right]}{(aw^2+bw+c)^m}dw=C\left(\frac{1}{5}\right)=
$$
\begin{equation}
=\left(\frac{2}{11+5 \sqrt{5}}\right)^{11/30} \frac{\Gamma \left(\frac{11}{30}\right) \Gamma \left(\frac{5}{6}\right)}{\Gamma \left(\frac{6}{5}\right)} \, _2F_1\left(\frac{1}{6},\frac{11}{30};\frac{6}{5};-\frac{123}{2}+\frac{55 \sqrt{5}}{2}\right).
\end{equation}

\section{Applications}

\subsection{\textbf{The study of $y(x)=F(x)=R\left(e^{-\pi\sqrt{k_i(x)}}\right)$ function}}

Assume that $y(x)=F(x)=R\left(e^{-\pi\sqrt{k_i(x)}}\right)$, then from Theorem's 1,5 we have 
$m_{G_0}(x)=k_i(x)$ and 
$$
5\int^{y(x)}_{0}\frac{dt}{t\sqrt[6]{t^{-5}-11-t^5}}=3\sqrt[3]{2x}\cdot {}_2F_1\left[\frac{1}{6},\frac{1}{3};\frac{7}{6};x^2\right]=\int^{x}_{0}\frac{dt}{G_0(y(t))}.
$$
Differentiating the above equation we get
\begin{equation}
h'_i(x)=\frac{1}{G_0(y(x))}=\frac{2^{1/3}}{x^{2/3}(1-x^2)^{1/3}}.
\end{equation}
Hence
\begin{equation}
y_i(x)=\sqrt{\frac{1-\sqrt{1-8G_0(x)^3}}{2}}
\end{equation}
and
$$
5\int^{x}_{0}\frac{G_0(t)}{t\sqrt[6]{t^{-5}-11-t^5}}dt=y_i(x).
$$ 
Hence the final equation for evaluating the $G-$function is
$$
\sqrt{\frac{1-\sqrt{1-8G_0(x)^3}}{2}}=5\int^{x}_{0}\frac{G_0(t)}{t\sqrt[6]{t^{-5}-11-t^5}}dt,
$$
which under differentiation becomes
$$
\frac{5}{x \sqrt[6]{-x^5+\frac{1}{x^5}-11}}=\frac{3 G_0(x) G'_0(x)}{\sqrt{\frac{1}{2}-4 G_0(x)^3} \sqrt{1-\sqrt{1-8 G_0(x)^3}}}.
$$
Solving the last differential equation, we get the following relation for $G_0(x)$:
$$
\frac{3\cdot 2^{2/3}\sqrt{G_0(x)}}{\sqrt[6]{1+\sqrt{1-8 G_0(x)^3}}}\cdot {}_2F_1\left[\frac{1}{6},\frac{1}{3};\frac{7}{6};\frac{1}{2}\left(1-\sqrt{1-8 G_0(x)^3}\right)\right]=
$$
\begin{equation}
=6x^{5/6} F_{Ap}\left[\frac{1}{6};\frac{1}{6},\frac{1}{6};\frac{7}{6};\frac{11-5 \sqrt{5}}{2}x^5,\frac{11+5 \sqrt{5}}{2}x^5\right]=F^{(-1)}_1(x).
\end{equation}
From (155) the partial evaluation of $G_0(x)$ follows.

\subsection{\bf The case of $G(x)=1$ function}

In case $G(x)=1$, then from relation (26) 
\begin{equation}
y(x)=F_1(x).
\end{equation} 
From $h(x)=x$, we have
\begin{equation}
F'_1(x)=5^{-1}F_1(x)\sqrt[6]{F_1(x)^{-5}-11-F_1(x)^5}
\end{equation}
and from (56)
\begin{equation}
Q'_G(x)=s'_i(x)=\frac{1}{\sigma(x)}=\frac{2^{1/3}}{\left(x\sqrt{1-x^2}\right)^{2/3}}.
\end{equation}
Hence
\begin{equation}
Q_G(x)=2^{1/3}\int\frac{dx}{\left(x\sqrt{1-x^2}\right)^{2/3}}=3\sqrt[3]{2x}\cdot{}_2F_1\left[\frac{1}{6},\frac{1}{3};\frac{7}{6};x^2\right]+c.
\end{equation}
The modular equation for $y_i(x)$ is $P_n(x)$ and
\begin{equation}
P_n(x)=5\int^{n^2F_1(x)}_{0}\frac{dt}{t\sqrt[6]{t^{-5}-11-t^5}}.
\end{equation}
Also
\begin{equation}
m_G(x)=b^{(-1)}_x=b^{(-1)}(x).
\end{equation}
Avoiding the inverse of $b_r$, which is a hypergeometric function (more precicely, Beta function) we can use the function $m_G(x)=m(x)$ defined by (see relation (23))
\begin{equation}
\pi\int^{+\infty}_{\sqrt{m(x)}}\eta\left(it/2\right)^4dt=x.
\end{equation}  
Hence
\begin{equation}
y(x)=F_1(x)=R\left(e^{-\pi\sqrt{m(x)}}\right)
\end{equation}
and
$$
m(x)=k_i\left(s(x)\right).
$$
Hence
\begin{equation}
Q\left(k_r\right)=\pi\int^{+\infty}_{\sqrt{r}}\eta\left(it/2\right)^4dt=b^{(-1)}_r=m(r).
\end{equation}
\begin{equation}
F_1^{(-1)}{'}(x)=\frac{5}{t\sqrt[6]{t^{-5}-11-t^5}}.
\end{equation}
$$
\int F_1(t)^{\nu}dt=\int \frac{x^{\nu}}{x\sqrt[6]{x^{-5}-11-x^5}}dx=
$$
\begin{equation}
=\frac{30x^{5/6+\nu}}{5+6\nu}F_{Ap}\left[\frac{1}{6}+\frac{\nu}{5},\frac{1}{6},\frac{1}{6},\frac{7}{6}+\frac{\nu}{5},\frac{11-5\sqrt{5}}{2}x^5,\frac{11+5\sqrt{5}}{2}x^5\right],
\end{equation}
where we have make the change of variable $x=F_1(t)$. Hence seting $\rho_1=\frac{11-5\sqrt{5}}{2}$, $\rho_2=\frac{11+5\sqrt{5}}{2}$ we get
\begin{equation}
\int F_1(t)^{\nu}dt=\frac{5F_1(t)^{5/6+\nu}}{5/6+\nu}F_{Ap}\left[\frac{1}{6}+\frac{\nu}{5},\frac{1}{6},\frac{1}{6},\frac{7}{6}+\frac{\nu}{5},\rho_1 F_1(t)^5,\rho_2 F_1(t)^5\right]+c.
\end{equation}
Hence if $G(t)=t^{\nu}$, then
$$
y\left(\frac{5F_1(c_A)^{5/6+\nu}}{5/6+\nu}F_{Ap}\left[\frac{1}{6}+\frac{\nu}{5},\frac{1}{6},\frac{1}{6},\frac{7}{6}+\frac{\nu}{5},\rho_1 F_1(c_A)^5,\rho_2 F_1(c_A)^5\right]\right)=F(A).\eqno{(167.1)}
$$
More general we have: If 
$$
G(t)=\sum^{\infty}_{n=1}\frac{G^{(n)}(0)}{n!}t^n,
$$
then
$$
\int G(F_1(t))dt=5F_1(c_A)^{5/6}\times
$$
$$
\times\sum^{\infty}_{n=0}\frac{G^{(n)}(0)}{n!}\frac{F_1(c_A)^{n}}{5/6+n}F_{Ap}\left[\frac{1}{6}+\frac{n}{5},\frac{1}{6},\frac{1}{6},\frac{7}{6}+\frac{n}{5},\rho_1 F_1(c_A)^5,\rho_2 F_1(c_A)^5\right]
$$
and
$$
y_i(F(A))=5F_1(c_A)^{5/6}\times
$$
$$
\sum^{\infty}_{n=0}\frac{G^{(n)}(0)}{n!}\frac{F_1(c_A)^{n}}{5/6+n}F_{Ap}\left[\frac{1}{6}+\frac{n}{5},\frac{1}{6},\frac{1}{6},\frac{7}{6}+\frac{n}{5},\rho_1 F_1(c_A)^5,\rho_2 F_1(c_A)^5\right].\eqno{(167.2)}
$$
Also 
$$
Q_{G}(A)=5F_1(c_A)^{5/6}\times
$$
$$
\times\sum^{\infty}_{n=0}\frac{G^{(n)}(0)}{n!}\frac{F_1(c_A)^{n}}{5/6+n}F_{Ap}\left[\frac{1}{6}+\frac{n}{5},\frac{1}{6},\frac{1}{6},\frac{7}{6}+\frac{n}{5},\rho_1 F_1(c_A)^5,\rho_2 F_1(c_A)^5\right].
$$
Hence
$$
h(A)=Q_G(c_i(A))=5F_1(A)^{5/6}\times
$$
$$
\times\sum^{\infty}_{n=0}\frac{G^{(n)}(0)}{n!}\frac{F_1(A)^{n}}{5/6+n}F_{Ap}\left[\frac{1}{6}+\frac{n}{5},\frac{1}{6},\frac{1}{6},\frac{7}{6}+\frac{n}{5},\rho_1 F_1(A)^5,\rho_2 F_1(A)^5\right].\eqno{(167.3)}
$$
This can be seen as the evaluation of Theorem 4 (evaluation with inverse integrals). However now we have infinite series expansion. The above equations can also help find $y(A)$ in case $G(A)$ is a polynomial.

\subsection{\bf The Case of $G(F(x))=x$ function}

From (77) we have
\begin{equation}
\sigma(x)=\frac{\left(x\sqrt{1-x^2}\right)^{2/3}}{\sqrt[3]{2}\cdot x}
\end{equation}
and
\begin{equation}
Q_G(x)=s_i(x)=\frac{3\sqrt[3]{x^4}}{2\sqrt[3]{4}}\cdot {}_2F_{1}\left[\frac{1}{3},\frac{2}{3};\frac{5}{3};x^2\right]+c.
\end{equation}
Hence
\begin{equation}
y\left(\frac{3\sqrt[3]{x^4}}{2\sqrt[3]{4}}\cdot {}_2F_{1}\left[\frac{1}{3},\frac{2}{3};\frac{5}{3};x^2\right]+c\right)=F(x).
\end{equation}
$$
h_1^{(-1)}(x)=G\left(y(x)\right)=G\left(F(k_{m_G})\right)=k(m_G(x))=Q_G^{(-1)}(x).
$$
Hence $h_1(x)=Q_G(x)$
$$
h'_1(x)=s'_i(x)=\frac{1}{\sigma(x)}=\frac{\sqrt[3]{2}\cdot x}{(x\sqrt{1-x^2})^{2/3}}.
$$

\subsection{\bf The Case of Jacobi Theta Functions}

For an extended version of this section one can see [15].\\
For $a,p$ positive rationals with $a<p$ and $|q|<1$ we define the forms
\begin{equation}
\left[a,p;q\right]:=\prod^{\infty}_{n=0}\left(1-q^{np+a}\right)\left(1-q^{np+p-a}\right)
\end{equation}
and
\begin{equation}
\theta(a,p;q):=\left[a,p;q\right]^{*}:=q^{C_0}\left[a,p;q\right],
\end{equation}
where $C_0:=p/12-a/2+a^2/(2p)$.\\

With the help of Jacobi triple product identity it can be shown that 
\begin{equation}
\vartheta\left(\frac{p}{2},\frac{p}{2}-a;q\right)=q^{C_0}\eta_1(p\tau)\left[a,p;q\right],
\end{equation} 
where
\begin{equation}
\vartheta\left(a,b;q\right):=\sum^{\infty}_{n=-\infty}(-1)^nq^{an^2+bn}\textrm{, }|q|<1,
\end{equation}
is a theta function and
\begin{equation}
\eta_1(\tau):=\prod^{\infty}_{n=1}\left(1-q^n\right)\textrm{, }q^{i\pi\tau}\textrm{, }\tau=i\sqrt{r}\textrm{, }r>0,
\end{equation}
is the Ramanujan-Dedekind eta function.\\

We also define
\begin{equation} 
Q_{\{a,p\}}(x):=\left[a,p;e^{-\pi\sqrt{k_i(x)}}\right]^{*}.
\end{equation}
Here $Q_{\{a,p\}}(x)$ is the $Q_G(x)$ function defined in Section 3 relation (52) above. In the case of Jacobi theta functions describes their algebraic part (conjecture).\\
\\

We will try to characterize these functions $Q_{\{a,p\}}(x)$. For this, assume that $P_{n}$ is the $n$-th modular equation of $\theta(a,p;q)$, then 
\begin{equation}
\theta(a,p;q^n)=P_{n}\left(\theta(a,p;q)\right).
\end{equation}
Also assume that our conjecture holds, then 
$$
Q_{\{a,p\}}\left(k_{n^2r}\right)=P_{n}\left(Q_{\{a,p\}}(k_r)\right).
$$
By inverting $k_r$, we get
$$
Q_{\{a,p\}}\left(k_{n^2k_i(x)}\right)=P_n\left(Q_{\{a,p\}}(x)\right).
$$
Setting 
\begin{equation}
S_n(x):=k_{n^2k_i(x)},
\end{equation}
we lead to the next\\
\\ 
\textbf{Theorem 24.}\\
If the $n$-th modular equation of $\theta(a,p;q)$ is that of (177),  then 
\begin{equation}
k_{n^2k_i(x)}=S_n(x)=Q_{\{a,p\}}^{(-1)}\left(P_{n}\left(Q_{\{a,p\}}(x)\right)\right),\textrm{ }n=2,3,4,...
\end{equation}
If one manages to solve equation (179) with respect to $Q_{\{a,p\}}(x)$ for given $a,p$, then  
\begin{equation}
\sum^{\infty}_{n=-\infty}(-1)^nq^{pn^2/2+(p-2a)n/2}=q^{-\frac{p}{12}+\frac{a}{2}-\frac{a^2}{2p}}\eta(q^p)Q_{\{a,p\}}(k_r), \forall r>0
\end{equation}
and $Q_{\{a,p\}}(x)$ will be a root of a minimal polynomial of degree $\nu=\nu(a,p,x)$.\\ 

Note that in case of rational $x\in(0,1)$ and $a,p$ rational with $0<a,p$, then the degree $\nu$ is independent of $x$ and the minimal polynomial of $Q_{\{a,p\}}(x)$ will have integer coefficients.\\
\\
\textbf{Example.}\\
The 2nd degree modular equation of $A(1,4;q)$ is
\begin{equation}
16 u^8+u^{16}v^8-v^{16}=0.
\end{equation} 
If we solve with respect to $v$ we get $v=P_2(u)$, where $v=A(1,4;q^2)$ and $u=A(1,4;q)$. Moreover
\begin{equation}
P_2(w)=\frac{\left(w^{16}+w^4\sqrt{64+w^{24}}\right)^{1/8}}{2^{1/8}}.
\end{equation}
It is $n=2$ then hold (see [9])
\begin{equation}
k_{4r}=\frac{1-\sqrt{1-k_r^2}}{1+\sqrt{1-k_r^2}}.
\end{equation}
Hence
\begin{equation}
S_2(x)=k_{4k_i(x)}=\frac{1-\sqrt{1-x^2}}{1+\sqrt{1-x^2}}.
\end{equation}
Finally we get from the relation (179) of Theorem 24:
\begin{equation}
\frac{\sqrt[8]{Q_{\{1,4\}}(x)^{16}+Q_{\{1,4\}}(x)^4\sqrt{Q_{\{1,4\}}(x)^{24}+64} }}{\sqrt[8]{2}}=Q_{\{1,4\}}\left(\frac{1-\sqrt{1-x^2}}{1+\sqrt{1-x^2}}\right),
\end{equation}
which has indeed a solution 
$$
Q_{\{1,4\}}(x)=\sqrt[12]{\frac{4(1-x^2)}{x}}.
$$  
\\
\textbf{Note.}\\
We note that function $m(q)=k_r^2$ is implemented in program Mathematica. However a useful expansion is
\begin{equation}
k_r=\sqrt{m(q)}=4q^{1/2}\exp\left(-4\sum^{\infty}_{n=1}q^n\sum_{d|n}\frac{(-1)^{d+n/d}}{d}\right),
\end{equation}        
where $q=e^{-\pi\sqrt{r}}$, $r>0$.\\

Continuing we denote
\begin{equation}
\theta(q):=\theta_{\{a,p\}}(q)=q^{p/12-a/2+a^2/(2p)}\frac{\vartheta\left(\frac{p}{2},\frac{p-2a}{2};q\right)}
{f\left(-q^p\right)}\textrm{, }q=e^{-\pi\sqrt{r}}.
\end{equation}
In the case of Jacobi theta functions we set 
\begin{equation}
Q(x):=Q_{\{a,p\}}(x)=\theta_{\{a,p\}}\left(e^{-\pi\sqrt{k_i(x)}}\right).
\end{equation}
But it holds $y(x)=F\left(k\left(m_G(x)\right)\right)$ and 
$m_G(x)=k_i\left(Q_i(x)\right)$, hence $y(x)=F\left(Q_i(x)\right)$, inverting
\begin{equation}
y_i(x)=\theta_{\{a,p\}}\circ k_i\circ F_i(x).
\end{equation}
Hence we get the following\\ 
\\
\textbf{Theorem 25.}\\
If $q=e^{-\pi\sqrt{r}}$, $r>0$, then
\begin{equation}
y\left(\theta(q)\right)=R\left(q\right)\textrm{ and }m_G^{(-1)}(r)=\theta(q).
\end{equation}
\\
\textbf{Theorem 26.}\\ 
The $n-$th modular equation of $\theta(q)=\theta_{\{a,p\}}(q)$ is
\begin{equation}
P_n(x)=Q_{\{a,p\}}\left(k\left(n^2k_i\left(Q_{\{a,p\}}^{(-1)}(x)\right)\right)\right).
\end{equation}
Also $m_G(x)=k_i\left(Q_{\{a,p\}}^{(-1)}(x)\right)$, $m_G^{(-1)}(n)=P_{\sqrt{n}}\left(Q_{\{a,p\}}(1)\right)$, $n>0$
and
\begin{equation}
y(x)=R\left(e^{-\pi\sqrt{k_i\left(Q^{(-1)}_{\{a,p\}}(x)\right)}}\right).
\end{equation}
\\
\textbf{Theorem 27.}\\
If $q=e^{-\pi\sqrt{r}}$, $r>0$, then
\begin{equation}
\frac{d\theta(q)}{dr}=\frac{1}{\phi(r)}.
\end{equation}
\\
\textbf{Proof.}\\
From $m_G^{(-1)}(A)=\theta(q)$, $q=e^{-\pi\sqrt{A}}$ we get 
\begin{equation}
5\int^{R(q)}_{0}\frac{G(t)dt}{t\sqrt[6]{t^{-5}-11-t^5}}=\theta(q).
\end{equation}
After derivating the above relation and using (5), we get 
\begin{equation}
G\left(R(q)\right)q^{-5/6}f(-q)^4=\theta'(q).
\end{equation} 
Using (14), we get the result.\\ 
\\
\textbf{Notes.}\\
Assuming the above we have
$$
\theta'(q)=q^{-1}\eta(z)^4G\left(R(q)\right),
$$
where 
$$
\eta(z):=q^{1/24}f(-q)\textrm{, }q=e^{2\pi i z}\textrm{, }Im(z)>0,
$$
is the Dedekind eta function (see also relation (1)).\\
\\
\textbf{Conjecture 1.}\\
In the case $Q_G(x)=Q_{\{a,p\}}(x)$ we have that $G(R(q))$ is root of polynomial with integer coefficients.\\
\\
\textbf{Definition 4.}\\
We call theta function of the $G-$transformation of Theorem 1 the function $m_G^{(-1)}(A)$, where $q=e^{-\pi\sqrt{A}}$, $A>0$. By this way the definition of the theta functions is generalized and related with the $G-$transform.\\

Therefore for the function $m_G(A)$ holds
\begin{equation}
\theta\left(e^{-\pi\sqrt{m_G(A)}}\right)=A
\end{equation}
and
\begin{equation}
R\left(e^{-\pi\sqrt{m_G(A)}}\right)=y(A).
\end{equation}
\\
\textbf{Theorem 28.}\\For the modularity of $m^{(-1)}_G(A)$ we have the next relation
\begin{equation}
m^{(-1)}_{G}\left(\frac{1}{A}\right)=Q_G\left(\sqrt{1-Q^{(-1)}_G\left(m_G^{(-1)}(A)\right)^2}\right).
\end{equation}
\\
\textbf{Proof.}\\
From relation (53) and the idenity $k_{1/r}=k'_r$, we get the result.\\
\\
\textbf{Example.}\\
Suppose $a=1$, $p=4$, then
\begin{equation}
\sum^{\infty}_{n=-\infty}(-1)^nq^{2n^2+n}=q^{1/24} \eta(4\tau)Q_{\{1,4\}}(k_r).
\end{equation}
Then $Q_{\{1,4\}}(x)$ will be 
$$
Q(x)=\sqrt[6]{2}\sqrt[12]{\frac{1-x^2}{x}}.
$$
For a certain $G$ we have from (62) and (64.1): 
$$
\sigma(x)=-6\cdot2^{5/6}x^{13/12}(1-x^2)^{11/12}(1+x^2)^{-1}.
$$  
Hence
$$
y(x)=R\left(e^{-\pi\sqrt{k_i\left(\frac{1}{8}\left(-x^{12}+\sqrt{64+x^{24}}\right)\right)}}\right).
$$

\[
\]

\centerline{\bf References}\vskip .2in

\noindent

[1]: M. Abramowitz and I.A. Stegun. 'Handbook of Mathematical Functions'. Dover Publications, New York.(1972)\\

[2]: G.E. Andrews. 'Number Theory'. Dover Publications, New York.(1994)\\
 
[3]: G.E. Andrews, Amer. Math. Monthly, 86, 89-108.(1979)\\

[4]: J.V. Armitage, W.F. Eberlein. 'Elliptic Functions'. Cambridge University Press.(2006)\\

[5]: N.D. Bagis and M.L. Glasser. 'Integrals related with Rogers Ramanujan continued fraction and $q$-products'. arXiv:0904.1641v1 [math.NT].(2009)\\
 
[6]: B.C. Berndt. 'Ramanujan`s Notebooks Part I'. Springer Verlag, New York.(1985)\\

[7]: B.C. Berndt. 'Ramanujan`s Notebooks Part II'. Springer Verlag, New York.(1989)\\

[8]: B.C. Berndt. 'Ramanujan`s Notebooks Part III'. Springer Verlag, New York.(1991)\\

[9]: B.C. Berndt. 'Ramanujan`s Notebooks Part V'. Springer Verlag, New York.(1998)\\

[10]: I.S. Gradshteyn and I.M. Ryzhik. 'Table of Integrals, Series and Products'. Academic Press.(1980)\\  

[11]: N.N. Lebedev. 'Special Functions and their Applications'. Dover Pub. New York.(1972)\\

[12]: E.T. Whittaker and G.N. Watson. 'A course on Modern Analysis'. Cambridge U.P.(1927)\\

[13]: N.D. Bagis. 'Generalized Elliptic Integrals and Applications'.\\arXiv:1304.2315v2 [math.GM].(2013)\\ 

[14]: Bruce C. Berndt, Heng Huat Chan, Sen-Shan Huang, Soon-Yi Kang, Jaebum Sohn and Seung Hwan Son. 'The Rogers-Ramanujan Continued Fraction'. J. Comput. Appl. Math., 105 (1999), 9-24.\\ 

[15]: N.D. Bagis. 'On the Complete Evaluation of Theta Functions'. arXiv:1503.01141v4 [math.GM] 10 Mar 2021.\\ 

\end{document}